\documentclass[11pt,a4paper]{article}
\usepackage[latin1]{inputenc}
\usepackage{ amsthm,amsmath,amsfonts,amssymb}
\usepackage{epsfig}
\usepackage{colordvi,helvet,multicol}
\usepackage{epic}
\usepackage{mathptmx}
\usepackage{times}
\usepackage[T1]{fontenc}
\usepackage{verbatim}
\usepackage{graphics}
\usepackage{makeidx}
\usepackage{fancyhdr}
\makeindex
\makeatletter
\long\def\@caption#1[#2]#3{%        Abbildungstext auf Fussnotengroesse
  \par
  \addcontentsline{\csname ext@#1\endcsname}{#1}%
    {\protect\numberline{\csname the#1\endcsname}{\ignorespaces #2}}%
  \begingroup
    \@parboxrestore
    \if@minipage
      \@setminipage
    \fi
   \normalsize
    \footnotesize
    \@makecaption{\csname fnum@#1\endcsname}{\ignorespaces #3}\par
  \endgroup}
\makeatother
 \theoremstyle{plain}
\newcommand{\ee}{\end{equation}}

\newtheorem{definition}{Definition}[section]
\newtheorem{theorem}[definition]{Theorem}
\newtheorem{lemma}[definition]{Lemma}

\newtheorem{proposition}[definition]{Proposition}
\newtheorem{corollary}[definition]{Corollary}

\newcommand\CC{{\Bbb C}}

\newcommand\NN{{\Bbb N}}

\begin{document}

\title{ On the Amenability of   Compact and Discrete  Hypergroup Algebras
}
\author{Ahmadreza Azimifard
\footnote{Mathematics Department, 
Stony Brook University, Stony Brook NY, 11794-3651, USA.\hspace{1.cm} azimifard@math.sunysb.edu }
}

\date{}
 \maketitle
\begin{abstract}
Let $K$ be a   commutative compact  hypergroup and $L^1(K)$ the hypergroup
algebra.
We show that $L^1(K)$ is amenable if and only if $\pi_K$, the Plancherel
weight on  the dual space $\widehat{K}$, 
 is bounded. Furthermore, we  show that if $K$ is an infinite discrete hypergroup and there exists 
   $\alpha\in \widehat{K}$ which vanishes at infinity, then $L^1(K)$ is not amenable. In particular,  $L^1(K)$ fails to be  even $\alpha$-left amenable if   $\pi_K(\{\alpha\})=0$. 
  \end{abstract}

 {\bf{Introduction.}}
 			Let $K$ be a     commutative compact hypergroup, $\widehat{K}$ its dual
			space, and $L^1(K)$ the hypergroup algebra. More recently in \cite{A.S.N},  among other things,
  			we  showed    that when $K$ is a 	hypergroup of  conjugacy classes of a non-abelian
 			compact connected Lie  group $L^1(K)$,  in contrast to the   group case,   is not amenable.
 			The proof of this theorem,  which is mainly based on the structure  of underlying  
 			group, follows from the fact that   the Plancherel weight on $\widehat{K}$
			 tends to infinity and consequently  the approximate diagonal for $L^1(K)$ is not bounded.
  			In this paper, we  show that the statement  remains valid for  general 
  			 commutative  compact hypergroups. More 	precisely, we show that  
  			   $L^1(K)$ is amenable if and only if   the Plancherel weight on $\widehat{K}$ is bounded.  
 			And, similar to the group case \cite{Zhang}, we also show that closed ideals of $L^1(K)$  possess
 			 approximate identities. In addition, we generalize our recent results on polynomial 
 			 hypergroups \cite{azimifard.Math.Z} to discrete hypergroups. If  $K$ is      
 			 a (infinite) discrete hypergroup and  $\alpha\in \widehat{K}$ which vanishes at infinity, then   
  			  $L^1(K)$ is not amenable. Indeed,
 			 we show that if $\pi_K(\{\alpha\})=0$, then 
   			$L^1(K)$ is not even  $\alpha$-left amenable, and $L^1(K)$  fails to be amenable when  $\pi_K(\{\alpha\})>0$. 
  			Observer that  in the latter case  $L^1(K)$ might be $\alpha$-left amenable; see \cite{azimifard.Math.Z}.

    {\bf{Preliminaries.}}
    			Let $(K, p, \sim )$ denote a
                  locally compact commutative hypergroup with Jewett's axioms \cite{Jew75}, where
                  $p:K\times K\rightarrow M^1(K)$, $(x,y)\mapsto p(x,y)$,
                  and $\sim:K\rightarrow K$,  $x\mapsto \tilde{x}$,
                  specify  the convolution and involution on $K$ and $p{(x,y)}=p{(y,x)}$ for
                  every $x, y\in K$. Here
                  $M^1(K)$  stands for  the set of all probability measures on  $K$.

            Let $C_c(K)$  be the
    space of all  continuous functions on $K$ with the uniform norm
    $\| \cdot  \|_\infty$.  The translation of $f\in C_c(K)$ at  the point
            $x\in K$, $T_xf$, is defined by $T_xf(y)=\int_K  f(t)dp{(x,y)}(t)$,
            for every $y\in K$.
  %  The character space of $K$ is denoted by
  %   $$ \mathfrak{X}^b(K):=\left\{
  %                                   \alpha\in C^b(K): \alpha(e)=1, \;\omega(x,y)
  %                                  (\alpha)=\alpha(x)\alpha(y), \; \forall\;x,y\in K
  %                       \right\}.
  %  $$
%
%
    Let $(L^1(K), \|\cdot\|_1)$  denote the  usual  Banach $*$-algebra
    of integrable functions on $K$ with respect to its Haar measure $m$,
    where the convolution and involution of $f,g\in
    L^1(K)$ are given by  $ f*g(x)=\int_K f(y)T_{\tilde{y}}g(x)dm(y)$ ($m$-a.e.)
    and
    $f^\ast(x)=\overline{f(\tilde x)}$
    respectively. 
    If $K$ is discrete, then $L^1(K)$ has an identity element;
              otherwise $L^1(K)$ has a bounded approximate identity, i.e.
              there exists a bounded net $\{e_i\}_i$ of functions in $L^1(K)$, $\|e_i\|_1\leq M$,   $M>0$,
              such that $\|f \ast e_i-f\|_1\rightarrow 0$ as
              $i\rightarrow \infty$.
     The dual of  $L^1(K)$ can be identified with
              the usual Banach space
              $L^\infty(K)$, and its structure space is  homeomorphic to the
              character
              space of $K$, i.e.
    \[ \mathfrak{X}^b(K):=\left\{
                                     \alpha\in C^b(K): \alpha(e)=1, \;p(x,y)
                                    (\alpha)=\alpha(x)\alpha(y), \; \forall\;x,y\in
                                    K\right\}\]
         equipped with the compact-open topology.
        $\mathcal{X}^b(K)$ is a locally compact Hausdorff space.
        Let  $\widehat{K} $ denote the set of all  hermitian
        characters $\alpha$ in $\mathcal{X}^b(K)$,
        i.e. $\alpha(\tilde{x})=\overline{\alpha(x)}$ for every $x\in
        K$, with   a  Plancherel measure $\pi_K$.
%      The Hilbert spaces  $L^2(K)$ and $L^2(\widehat{K})$ are isometric.
     Observe  that   $\widehat{K}$ in general may not have the dual hypergroup structure  and  a
     proper inclusion in 
     $\mbox{supp}({\pi_K})\subseteq\widehat{K}\subseteq  \mathcal{X}^b(K)$   is possible. If $K$ is
      compact, then the dual space is unique and it is dense in $C(K)$ (see \cite{BloHey94, Jew75}).
     % If $K$ is compact, then the linear span of $\mbox{supp}({\pi})$ $(=\mathcal{X}^b(K))$ 
     %is dense in $L^2(K)$ \cite{Jew75}.    

           The Fourier-Stieltjes transform of $\mu\in M(K)$, $\widehat{\mu}\in C^b(\widehat{K})$,
       is  given by  $\widehat{\mu}(\alpha):=\int_K \overline{\alpha(x)}d\mu(x)$. 
       Its restriction to   $L^1(K)$   is called the Fourier transform.
       We have  $\widehat{f}\in C_0(\widehat{K})$, for $f\in
       L^1(K)$, and the map
$\alpha\rightarrow I(\alpha):=\mbox{ker}(\varphi_\alpha)$
 is a bijection of $\widehat{K}$ onto the
space of all maximal ideals of $L^1(K)$, where $\mbox{ker}(\varphi_\alpha)$
 denotes the kernel of the homomorphism
  $\varphi_\alpha(f)=\widehat{f}(\alpha)$  on $L^1(K)$; see \cite{BonDun73}.  
  
 Let $X$ be a   Banach $L^1(K)$-bimodule and $\alpha\in
               \widehat{K}$. In a canonical way
                the dual space $X^\ast$ is a Banach
                $L^1(K)$-bimodule. The module   $X$ is called a $\alpha$-left
               $L^1(K)$-module if the left module multiplication is given by
               $f\cdot x=\widehat{f}(\alpha)x$, for  every
                $f\in L^1(K)$ and $x\in X$. In this case,  $X^\ast$ turns out  to be a
                $\alpha$-right $L^1(K)$-bimodule,
                i.e.  $\varphi\cdot f=\widehat{f}(\alpha)\varphi$,  for
                every $f\in L^1(K)$ and $\varphi\in X^\ast$.
A continuous  linear map $D:L^1(K)\rightarrow X^\ast$
                is called a derivation if $D(f\ast g)=D(f)\cdot g+f\cdot D(g)$,
                for every  $f,g\in L^1(K)$, and  an inner derivation
                if $D(f)=f\cdot \varphi-  \varphi\cdot f$,
                for some $\varphi\in X^\ast$.   The algebra $L^1(K)$ is called
                $\alpha$-left amenable if for every $\alpha$-left
               $L^1(K)$-module $X$,    every
               continuous derivation $D:L^1(K)\rightarrow X^\ast$
               is inner; and, if the latter holds for every
               Banach $L^1(K)$-bimodule $X$,
               then $L^1(K)$ is called amenable.\\

				Let $K'=K\times K$ denote the hypergroup of  cartesian product of  
				$K$ with itself.  It is straightforward to show that  $L^1(K')\cong L^1(K)\otimes_pL^1(K)$
                        ($\otimes_p$ denotes the projective tensor product) and 
              with  the actions 
			$f\cdot(g\otimes h)=(f\ast g)\otimes h$ and $(g\otimes h)\cdot f=g\otimes (h\ast f)$ 
			the Banach algebra  $L^1(K')$ becomes a $L^1(K)$-bimodule. We observe  that 
  			the map $\phi:\mathcal{X}^b(K)\times \mathcal{X}^b(K)\rightarrow \mathcal{X}^b(K')$ defined 
  			 by $(\alpha,\beta)\rightarrow \alpha\otimes \beta $
  			 is a  homeomorphism (see \cite{BonDun73}). 
 			 As shown in \cite{Johnson},  $L^1(K)$ is amenable if it admits a bounded approximate diagonal,
 				i.e. a bouned net $\{M_i\}_i\subset L^1(K)\otimes_p L^1(K)$  which satisfies
		 \begin{equation}\notag
		 \pi(M_i)\cdot f, f\cdot \pi(M_i)\rightarrow f \mbox{ and } {f\cdot M_i-M_i\cdot f\rightarrow f}
		 \end{equation}
		for any $f\in L^1(K)$, where $\pi:L^1(K)\otimes_p L^1(K) \rightarrow
		L^1(K)$ is the convolution map. The amenability
		of $L^1(K)$ is also equivalent to the existence of a virtual diagonal,
		i.e. an element $M\in (L^1(K)\otimes_pL^1(K))^{\ast\ast}$ such that
	\begin{equation}\notag
		f\cdot M=M\cdot f\hspace{.5cm} f\pi^{\ast\ast}(M)=\pi^{\ast\ast}(M) f=f
	\end{equation}
	for any $f\in L^1(K)$, where the module actions of $L^1(K)$ on $(L^1(K)\otimes_pL^1(K))^{\ast\ast}$ and $L^1(K)^{\ast\ast}$
	are the second adjoints of the module actions of $L^1(K)$ on $L^1(K)\otimes_pL^1(K)$ and $L^1(K)$, respectively, and
	$\pi^{\ast\ast}$ is the second adjoint of $\pi$. We also define   $\pi_1,\pi_2:L^1(K)\rightarrow L^1(K')$ by 
	$\pi_1(f)(x,y)=f(x)\delta_e(y)$ and $\pi_2(f)=f(y)\delta_e(x)$, respectively, when $K$ is discrete. One can easily verify that 
	 the $\pi_i$ maps   are isometric and 
	$\pi_i(f\ast g)=\pi_i(f)\ast \pi_i(g)$ for every $f,g\in L^1(K)$.\\ 

As already mentioned, in this paper %in  continuing our work   in \cite{azimifard.Math.Z, A.S.N}, 
  we deal with the amenability problem 
   of compact and discrete hypergroup algebras. The results are organized as follows. 
     We first     show that a compact hypergroup algebra 
     $L^1(K)$ is amenable if and only if the Plancherel weight 
     $\pi_K$ on $\widehat{K}$ is bounded (Theorem  \ref{main.2}).
     %   $L^1(K)$ is weakly amenable(Proposition  \ref{q.2}).    
    Moreover, we  show that every closed ideal of $L^1(K)$ has an approximate identity 
     (Theorem \ref{closed.1}). We then discuss  amenability of non-compact discrete hypergroup 
     algebras. Let $K$ be a discrete hypergroup and $\alpha\in \widehat{K}$.
      If $\alpha$ vanishes at infinity, then $L^1(K)$ is not amenable; in the case of  $\pi_K(\{\alpha\})=0$, particularly, 
       the algebra  $L^1(K)$ is not even $\alpha$-left amenable (Theorem \ref{theorem.2}). Using our theorems, we finally  examine  the amenability of hypergroup algebras  of various  compact and discrete hypergroups.  
       
       I would like to thank Dr. Nico Spronk for  his comment on the early draft of this paper. 

\section{Amenability of Compact Hypergroup Algebras}
 
%   As we have  already shown in \cite{azimifard.Math.Z}, a compact hypergroup algebra 
%    $L^1(K)$ is $\alpha$-left amenable in every $\alpha\in \widehat{K}$, however
As it is already shown in  \cite{A.S.N},  if $K$ is a hypergroup of  conjugacy classes of a compact connected Lie 
 group, then $L^1(K)$
    is amenable if and only if the dimension  of   irreducible 
 unitary representations of the group  is bounded. In the following theorem 
 we show that   the  statement remains valid in general.

 \begin{theorem}\label{main.2}
                            \emph{  Let   $K$ be   a compact hypergroup. Then $L^1(K)$ is amenable if and only
                            if the Plancherel weights on $\widehat{K}$  is bounded, i.e.,  there exists a $c>0$ such that
                               $\pi_K(\{\alpha\}) <c$ for all  $\alpha\in \widehat{K}$.}
            \end{theorem}
            
            Before proceeding to the proof of this theorem, let us first discuss the existence of and pertinent topics to 
            the  approximate diagonals 
            for compact hypergroup algebras. 
            
           We observe that  since the convolution map $(x,y)\rightarrow p(x,y)$, $K'\rightarrow M^1(K)$, is
            continuous ($M^1(K)$ is considered with the $\mbox{weak}^\ast$ topology), 
            a hypergroup algebra $L^1(K)$ is  $\mbox{weak}^\ast$ dense in $M(K)$, and
            the convolution map $\pi:L^1(K')\rightarrow L^1(K)$ has a $\mbox{weak}^\ast$
             extension $\tilde{\pi}:M(K')\rightarrow M(K)$ which is defined by 
             \begin{equation}\label{homo}\notag
             \int_Kf(x)d\tilde{\pi}(\mu)(x)= \int_{K'} T_xf(y)d\mu(x,y)\hspace{1.cm}f\in C(K).
             \end{equation}
             Obviously  we have $\tilde{\pi}(\mu\otimes \nu)=\mu\ast \nu$,   $\mu,\nu\in M(K)$, 
              and if for  a $f\in C(K)$ we 
             let $g(x,y)=T_xf(y)$, then $g\in C(K')$ and 
             \begin{align}\notag
             \tilde{\pi}(\mu\ast \nu)(f)&=\int_{K'}T_xf(y)d\mu\ast \nu(x,y)\\ \notag
             &=\int_{K'}\int_{K'}T_{(x_1,x_2)}g(y_1,y_2)d\mu(x_1,x_2)d\nu(y_1,y_2) \\ \notag
             &=\int_{K'}\int_{K'} T_{y_1}(T_{x_2}  T_{x_1}f)(y_2)d\nu(y_1,y_2)d\mu(x_1,x_2)\\ \label{homomorphism}
             &=\tilde{\pi}(\mu)\ast\tilde{\pi}(\nu)(f).
             \end{align}
               Hence  $\tilde{\pi}$ is a  homomorphism.

\begin{lemma}\label{theorem.1}
       \emph{ Let  $\{e_n\}$ be  a bounded
              approximate identity for $L^1(K)$, where
              $e_n=\sum_{m=0}^\infty a_m^n\alpha_m$ such that $a_m^n=0$ except
               for finitely many $m$.
Then
\begin{enumerate}
 \item[(i)]   $a_m^n\rightarrow \frac{1}{\|\alpha_m\|_2^2}$, and
            \item[(ii)] $M_n=\sum_{m=0}^\infty \left(a_m^n\right)^2\alpha_m\otimes \alpha_m$
                         is an  approximate diagonal for $L^1(K)$.
          %  \item[(iii)] if $L^1(K)$ is amenable,
          %               then there exists a measure $\mu\in M(K')$
          %               such that $\widehat{\mu}(\alpha,\beta)=\delta_\alpha(\beta).$
            \end{enumerate}
            }
\end{lemma}
            \begin{proof}
                            Let $\{U'_n\}$ be a family of neighborhoods of
                            the identity element $e$. Then the sequence
                            $\{e_n\}=\{\frac{1}{m(U'_n)}\chi_{U'_n}\}$ is a bounded
                            approximate identity for $L^1(K)$. Since the linear
                            span of $\widehat{K}$ is dense in $L^1(K)$, we may choose
                            $e_n=\sum_{m=0}^\infty a_m^n\alpha_m$, where  $a_m^n=0$
                            except for finitely many $m$. Therefore,
        \begin{equation}\notag
            \|\alpha_i\|_1\left|1-\widehat{e_n}(\alpha_i)\right|=
            \left\|\alpha_i-\widehat{e_n}(\alpha_i)\right\|_1
            =\left\|\alpha_i-e_n\ast \alpha_i\right\|_1\rightarrow 0\hspace{.5cm}(n\rightarrow \infty),
        \end{equation}
                     which implies that
                      $\|\alpha_i\|_1\left| 1-a_i^n\|\alpha_i\|_2^2\right|\rightarrow 0$,
                      consequently  $a_i^n\rightarrow \frac{1}{\|\alpha_i\|_2^2}$
                      as $n\rightarrow \infty.$

We now show that  $M_n=\sum_{m=0}^\infty\left(a_m^n\right)^2\alpha_m\otimes \alpha_m$ is an approximate diagonal for $L^1(K)$. Since
               \begin{equation}\notag
                            \pi(M_n)=\sum_{m=0}^n\left(a_m^n\right)^2\alpha_m\ast \alpha_m
                            =\sum_{m=0}^\infty (a_m^n)^2\|\alpha_m\|_2^2\alpha_m=e_n\ast e_n
            \end{equation}
        which is also a bounded approximate identity for $L^1(K)$ and
        \begin{align}\notag
                        \alpha_k\cdot M_n&=
                        \sum_{m=0}^\infty\left(a_m^n\right)^2\alpha_k\ast\alpha_m\otimes\alpha_m
                        =
            \sum_{m=0}^\infty\delta_k(m)\left(a_m^n\right)^2\alpha_m\otimes \alpha_m\\ \notag
            &=\sum_{m=0}^\infty\left(a_m^n\right)^2 \alpha_m\otimes \alpha_m\ast \alpha_k
            =M_n\cdot \alpha_k,
            \end{align}
            $\{M_n\}$ is an approximate diagonal
            for $L^1(K)$. Therefore, if   $\{M_n\}_n$ is bounded,
            then $L^1(K)$ is amenable \cite{Johnson}.
\end{proof}

We  now  use the idea in  the proof of  \cite[Theorem 1.6]{A.S.N} to establish the following lemma in its analogy.  

\begin{lemma}\label{Amenability.Key}
\emph{Let $K$ be a compact hypergroup and $\{M_n\}$ as in Lemma \ref{theorem.1}. Then the following statements are equivalent: 
 \begin{itemize}
 \item[(i)]$L^1(K)$ is   amenable.
 \item[(ii)]$\{M_n\}_n$ is bounded.
 \item[(iii)] There exists a measure $\mu\in M(K')$
                         such that $\widehat{\mu}(\alpha,\beta)=\delta_\alpha(\beta)$, $\tilde{\pi}(\mu)=\delta_e$, 
                         and $(f\otimes\delta_e)\ast \mu=\mu\ast(\delta_e\otimes f)$ for any $f\in L^1(K)$. 
 \end{itemize}}
\end{lemma}
\begin{proof}
$(i)\rightarrow (ii)$. In this case $L^1(K)$ admits a  bounded approximate diagonal, say
              $\{M'_k\}$. Let us assume that $M$ is the virtual diagonal and  $M'_k\overset{w^\ast}{\rightarrow} M$
             in $L^1(K')^{\ast\ast}$.
             Suppose  $\{e_n\}$ to  be as above and $F_n:=\{\alpha_m; a_m^n\not=0\}$.
             Then $F_n\otimes F_n$ is a finite dimensional ideal in
             $L^1(K)\otimes L^1(K)$ which contains $e_n\otimes e_n$.
             Then $\{e_n\otimes e_n\ast M'_k\}$ is a bounded  net in
             $A_n\otimes A_n,$ $A_n=\langle F_n\rangle$, with a limit point $N_n$. Write
             $N_n= {\sum}_{\alpha_i,\alpha_j\in F_n}c^n_{ij}\alpha_i\otimes \alpha_j$.
             For every $\alpha_m\in F_n$, since $M_k\cdot \alpha_m=\alpha_m\cdot M_k$ for every $k$,  we have 
             $\alpha_m\cdot N_n=N_n\cdot \alpha_m$. Therefore   
             %
             %
            %\begin{align}\notag
             %\[       \alpha_m\cdot N_n=
                              %  e_n\otimes e_n\ast\left[\underset{k\rightarrow \infty }{\lim}
                              %  \alpha_m\cdot M'_k\right]
                              %  =e_n\otimes e_n\ast \left[\underset{k\rightarrow \infty }{\lim}
                              %  M'_k\cdot \alpha_m
                              %                      \right]
                              %                        = \underset{k\rightarrow \infty }
                              %  {\lim}\left[e_n\otimes e_n\ast M'_k\right]\cdot \alpha_m
             %                   =N_n\cdot \alpha_m
             %                   \]
           % \end{align}
%  Therefore,
%
            \[
                    \sum_{\alpha_i,\alpha_j\in F_n}c^n_{ij}\delta_i(m)
                    \|\alpha_i\|_2^2\alpha_i\otimes\alpha_j=
                    \sum_{\alpha_i,\alpha_j\in F_n}c^n_{ij}\delta_j(m)
                    \|\alpha_j\|_2^2 \alpha_i\otimes\alpha_j
            \]
    which implies
    $\sum_j c^n_{mj}\|\alpha_m\|_2^2 \alpha_m\otimes \alpha_j
    =\sum_{i}c^n_{im}\|\alpha_m\|_2^2 \alpha_i\otimes \alpha_m$. Hence,
     from the orthogonality   of characters it follows  that
    $c^n_{mj}=0$ if $m\not=j$, so $N_n=\sum c_{ii}^n\alpha_i\otimes \alpha_i$. We have 
                \[
                        \pi(N_n)
                        =\pi(e_n\otimes e_n)\ast \underset{k\rightarrow \infty}{\lim}{\pi(M_k)}
                        =\pi(e_n\otimes e_n)=e_n\ast e_n,
                \]
         and in particular 
                \[
                        \sum_i c_{ii}^n\|\alpha_i\|_2^2\alpha_i
                        =\sum_i(a_i^2)^2\|\alpha_i\|_2^2\alpha_i,
                \]
            which yields  $c_{ii}^n=(a_i^n)^2$ for
            each $i$.  Hence   $M_n=N_n$ and  boundedness of $\{\|M_n\|_1\}$ follows from  
            $\|M_n\|_1=\|N_n\|_1\leq \|e_n\|_1^2\underset{k\rightarrow \infty}{\sup}
            \|M_k\|_1<\infty$.
            % provide the boundness of $\{\|M_n\|_1\}$.
           % Therefore, $\{M_n\}$ is the only approximate diagonal for $L^1(K')$ and if $\{\|M_n\|_1\}_n$ is bouned then $L^1(K)$ is amenable.

           $(ii)\rightarrow (iii)$. Since the algebra $L^1(K')$ can be canonically embedded in $M(K')$, it 
           follows from Banach-Alaoglu's theorem  that $\{M_n\}_n$  has a $\mbox{weak}^\ast$ limit point   $M\in M(K')$. 
%            Since
 %                    $M_n\in L^1(K')\subseteq M(K')\cong C(K')^\ast$,
                     % and $\{M_n/\|M_n\|_1\}$
                    %is a subset of the unit ball in $C(K')^\ast$, which is  $w^\ast$-compact
  %                   there exists  $\mu\in M(K')$
   %                 such that $M_n\overset{w^\ast}{\rightarrow} \mu$ when 
    %                $n\rightarrow \infty$; see \cite{HewRos79I}. 
    We have 
           \begin{align}\notag
                       \widehat{\mu}(\alpha_m,\alpha_{m'})&
                       =\underset{n\rightarrow \infty}{\lim}M_n(\alpha_m \otimes \alpha_{m'})
                       = \underset{n\rightarrow \infty}{\lim}\int_{K'} \sum_{i=0}^\infty
                       \left(a_i^n\right)^2\alpha_i(x)\alpha_i(y)
                       \overline{\alpha_m(x)} \overline{\alpha_{m'}(y)}dm(x)dm(y)\\ \notag
                       &=\underset{n\rightarrow
                       \infty}{\lim}\sum_{i=0}^\infty\left(a_i^n\right)^2\left(\int_K
                       |\alpha_i(x)|^2dm(x)\right)
                       \left(\int_K |\alpha_i(y)|^2dm(y)\right)
                       \delta_i(m)\delta_i(m')=\delta_m(m'). \notag
                       %&=\sum_{i=0}^\infty \underset{n\rightarrow \infty}
                       %{\lim}\left(a_i^n\right)^2 \| \alpha_i \|_2^2
                       %\delta_i(m) \delta_i(m')=\delta_m(m')\notag
         \end{align}
In that $\widehat{M(K')}\subseteq C^b(\widehat{K'})$, we now define the map
 $D:C^b(\widehat{K}\times \widehat{K})\rightarrow C^b(\widehat{K})$
by $D\mu(\alpha)=\widehat{\mu}(\alpha,\alpha)$. Obviously  for any $\nu\in M(K')$ we have 
$\widehat{\tilde{\pi}(\nu)}(\alpha)=D\widehat{\nu}(\alpha)$  and, in particular,  
\[ {\tilde{\pi}(\mu)}^{\hat{\;}}(\alpha)=D\widehat{\mu}(\alpha)=1=\widehat{\delta_e}(\alpha) \hspace{.5cm}(e\in K)\]
It follows from  the inverse of the Fourier transform \cite{BloHey94}
 that  $ {\tilde{\pi}(\mu)}=\delta_e$. We see, in addition,   that  
  if $f\in L^1(K)$  and $\alpha\in \widehat{K}$, ${(f\otimes \delta_e)}^{\hat{\;}} (\alpha,\beta)=\widehat{f}(\alpha)$ and
  $ {(\delta_e\otimes f)}^{\hat{\;}}(\alpha,\beta)=\widehat{f}(\beta)$. 
 Therefore $(f\otimes\delta_e)\ast \mu=\mu\ast(\delta_e\otimes f). $
 
 $(iii)\rightarrow (i).$ Let $\{e'_n\}_n$ be a bounded approximate identity in $L^1(K')$ and assume $M$ to be  
  a $\mbox{weak}^\ast$-limit point of $\{\mu\ast e'_n\}_n$ in $L^1(K')$. We shall show that $M$ is a virtual diagonal.
For any $f\in L^1(K)$ we have 
\[f\cdot M=\underset{n}{\lim}(f\otimes \delta_e)\ast \mu \ast e'_n=
\underset{n}{\lim}\mu\ast(\delta_e\otimes f)\ast e'_n=\underset{n}{\lim}\mu\ast e'_n\ast(\delta_e\otimes f)=M\cdot f.\]  
  And, if
   $E$ is  a $\mbox{weak}^\ast$-limit point of $\{\pi(e'_n)\}$, from 
  $\tilde{\pi}(\mu)=\delta_e$ and   (\ref{homomorphism})  it follows that   
  \[\pi^{\ast\ast}(M)=\underset{n}{\lim}\pi(\mu\ast e'_n)=
   \underset{n}{\lim}\tilde{\pi}(\mu)\ast\pi( e'_n)=\underset{n}{\lim}\pi( e'_n)=E. \]
  We obviously see that  $f\cdot E=E\cdot f=f$ for any $f\in L^1(K)$.   Therefore $M$ is a virtual diagonal.   
            \end{proof}

 %%%%%%%%%%%%%%%%%%%%%%%%%%%%%%%%%%%%%%%%%%%%%%%%%%%%%%%%%%%%%%%%%%%%%%%%%

           We now prove Theorem \ref{main.2} as follows: 
    \begin{proof}[{\bf{Proof of Theorem \ref{main.2}}}]
    First assume    that $L^1(K)$ is   amenable and    in contrary there exists a 
    sequence $\{\alpha_i\}_{i\in \NN}\subset \widehat{K}$ such that  $\pi_K(\{\alpha_i\})\rightarrow \infty$ 
    %$(\mbox{ obviously } \pi_K(\{\alpha_i\})>0)$
                          as $i\rightarrow \infty$. Obviously $ \pi_K(\{\alpha_i\})>0$ and 
                          the  functionals $F_{\alpha_i}:\widehat{K}\rightarrow \CC$  defined  by  
                           $F_{\alpha_i}(\beta)=\delta_{\alpha_i}(\beta)$ belong to 
                           $ L^1(\widehat{K})\cap L^2(\widehat{K})$. It is worth noting  that 
                           by the inverse of  Fourier transform we have 
\[\check{F_{\alpha_i}}(x)=\int_{\widehat{K}}F_{\alpha_i}(\beta)\beta(x)d\pi_K(\beta)
=\alpha_i(x)\pi_K(\{\alpha_i\}),\]
and  from  Plancherel's  theorem (see \cite{BloHey94}) we deduce   that 
$\pi(\{\alpha_i\})=\frac{1}{\|\alpha_i\|_2^2}>0$. By previous theorem there exists 
                         a $\mu\in M(K')$ such that
                         \[
                         1=\underset{i\rightarrow
                                    \infty}{\lim}\widehat{\mu}(\alpha_i,\alpha_i)
                                =\underset{i\rightarrow \infty}{\lim}\int_{K'}{\overline{\alpha_i(x)}}{\overline{\alpha_i(y)}}d\mu(x,y)
                                =\int_{K'}\underset{i\rightarrow
                                \infty}{\lim}\overline{\alpha_i(x)\alpha_i(y)}d\mu(x,y)=0,
                         \]
                         which is a contraction.

               To prove the converse of  the  theorem, let  $ {\sup}_{\alpha\in \widehat{K}}\pi_K(\{\alpha\}) <c$ for some $c>0$. 
              Since  $\{M_n\}$ is an approximate diagonal
                for $L^1(K')$ (Lemma  \ref{theorem.1}), by previous lemma it    
                suffices to show 
                that   $\{M_n\}$ is bounded.   For any  $f, g\in C(K)$ we have 
               \begin{align}
               %\underset {n\rightarrow \infty}{\lim}
               \left|\underset {n\rightarrow \infty}{\lim} M_n(f \otimes  g) \right|&=
                                  \left| \underset {n\rightarrow \infty}{\lim}\int_{K'}\sum_{i=0}^\infty
                                              (a_i^n)^2\alpha_i\otimes \alpha_i(x,y)\overline{f(x)}\overline{g(y)}dm(x)dm(y)\right|\\ \notag
               &\leq
               \sum_{i=0}^\infty \pi_K(\{\alpha_i\})^2  |\langle \overline{f}, \alpha_i\rangle ||\langle \overline{g}, \alpha_i\rangle|
               \leq c^2  \sum_{i=0}^\infty  |\langle \overline{f}, \alpha_i\rangle | \langle \overline{g}, \alpha_i\rangle| 
               \hspace{.5cm}( \mbox{Lemma} \ref{theorem.1})
               \\ \notag
               &
               \leq c^2 \sum_{i=0}^\infty  |\langle \overline{f}, \alpha_i\rangle |^2\cdot
               \sum_{i=0}^\infty  |\langle \overline{g}, \alpha_i\rangle |^2\leq
               c^2   \|f\|_2^2\|g\|_2^2\leq c^2\|f\|_\infty\|g\|_\infty. \notag
               \end{align}
                The latter   inequality  follows from  Plancherel and Cauchy-Schwartz's theorems.
                  Therefore   $L^1(K)$ is  amenable.  
\end{proof}

%%%%%%%%%%%%%%%%%%%%%%%%%%%%%%%%%%%%%%%%%%%%%%%%%%%%%%%%%%%%%%%%%%%%%%%%%%%%

Following \cite{Bade.dales.curtes} we say $L^1(K)$ is weakly amenable if every continuous  
derivation of $L^1(K)$ into $L^\infty(K)$ is zero. In contrast to the amenability of $L^1(K)$ we show that 
$L^1(K)$ is always weakly amenable when $K$ is compact. 

%Since the linear span of $\widehat{K}$ is dense in $L^p(K)(p=1,2)$, regardless of 
%asymptotic behaviour of the Plancherel weight on $\widehat{K}$, we expect the weak 
%amenability of $L^1(K)$ in general. By definition \cite{Bade.dales.curtes},  $L^1(K)$ is called to be weakly amenable
% if every bounded derivation of $L^1(K)$ into $L^\infty(K)$ is zero.
%We may expect the following simple proposition.

\begin{proposition}\label{q.2}
 \emph{Let  $K$ be a compact hypergroup. Then $L^1(K)$ is weakly amenable.
}
 \end{proposition}
%\vspace{-4mm}
\begin{proof}
Let  $D: L^1(K)\rightarrow L^\infty(K)$ be a continuous  derivation.
Due  to $\alpha\ast\alpha= {\|\alpha\|_2^2}\alpha$,
 for every  $\alpha\in \widehat{K}$,  we have
$D(\alpha)=\left( {2}/{\|\alpha\|_2^2}\right) \alpha\cdot D(\alpha)$. Here $"\cdot"$ stands for an arbitrary 
module action of $L^1(K)$ to $L^\infty(K)$. Hence 
 \begin{align}\notag
\alpha\cdot D(\alpha)&=
\left(2/\|\alpha\|_2^2\right)\left[ \alpha\cdot (\alpha\cdot D(\alpha))\right]\\[.1cm]
\notag
 &=\left(2/\|\alpha\|_2^2\right)\left[(\alpha\ast\alpha)\cdot D(\alpha)\right]\\[.1cm]
\notag &=2\alpha\cdot D(\alpha)\notag
\end{align}
which implies that    $D(\alpha)=0$. Since the linear span of $\widehat{K}$ is dense in $L^1(K)$, we obtain   $D=0$, as desired. 
\end{proof}

%%%%%%%%%%%%%%%%%%%%%%%%%%%%%%%%%%%%%%%%%%%%%%%%%%%%%%%%%%
As already mentioned  since $L^1(K)$, a compact hypergroup algebra,   is $\alpha$-left amenable in
 every $\alpha\in \widehat{K}$, the maximal ideals of $L^1(K)$ possess bounded approximate identities; 
 see \cite[1.2]{azimifard.Math.Z}. In the sequel, similar to the compact group case in \cite{Zhang}, we show that 
 closed ideals in $L^1(K)$ contain approximate identities.  

\begin{lemma}\label{c.1}
\emph{ Let  $J$  be a closed ideal of
$L^1(K)$   and
 $I_\alpha:=\underset{\beta\not=\alpha}\bigcap I(\beta)$.
Then
\begin{itemize}
\item[(i)]   $I_\alpha \simeq \CC \alpha$,  for every
$\alpha\in \widehat{K},$
\item[(ii)]$I_\alpha \subseteq J$ if and only if
$\widehat{f}(\alpha)\not=0$,  for some $f\in J$, and 
\item[(iii)] the
map $\alpha\mapsto I_\alpha$ is bijective from $\widehat{K}$
    onto  the set of all minimal ideals of $L^1(K)$.
\end{itemize}}
\end{lemma}

\begin{proof}
 (i)
 Let $\alpha\in \widehat{K}$. Obviously $I_\alpha \cap I(\alpha)= \{0\}$ and
$\alpha\in I_\alpha\cap \left(L^1(K)\setminus I(\alpha)\right)$. Let
$f$ be a non-zero element in $ I_\alpha$. Then 
$\lambda=\widehat{f}(\alpha)\not=0$  and  $\widehat{\lambda\cdot \alpha}(\beta)=\left(\lambda
\|\alpha\|_2^2\right)\delta_\alpha(\beta)$ which implies that
$f=\frac{\lambda}{\|\alpha\|_2^2}\cdot \alpha$. Hence
$I_\alpha\simeq\CC \alpha$, as desired.

(ii) Suppose  $f\in J$  with $\widehat{f}(\alpha)\not=0$. Since
$f\ast \alpha \in I_\alpha\cap J$, $f\ast
\alpha=\widehat{f}(\alpha)\alpha\not=0$,  and $I_\alpha \simeq \CC
\alpha$, we have $I_\alpha\subseteq I_\alpha\cap J $; thus
$I_\alpha\subseteq  J$.

(iii) Since    $ J\not=\{0\}$, there  exist
$f\in J$ and
 $\alpha\in \widehat{K}$ such that $\widehat{f}(\alpha)\not=0$.
 By (ii)  we have $I_\alpha\subseteq J$, consequently $J=I_\alpha$ as  $J$ is a minimal ideal.
 \end{proof}

%%%%%%%%%%%%%%%%%%%%%%%%
\begin{corollary}\label{c.4}
\emph{The proper closed ideals
of $L^1(K)$ are exact the family $\{I_P:\hspace{.2mm}P\subset
\widehat{K}\}$, where $I_P$ denotes  the closure of the linear span
of $P $ in $L^1(K)$.  Different closed   subsets of $\widehat{K}$
generate in this way  different closed  ideals.}
\end{corollary}

\begin{theorem}\label{closed.1}
\emph{Let $K$ be a compact    hypergroup.
 Then  every closed ideal
of $L^1(K)$ has an approximate identity.
  }
\end{theorem}
\begin{proof}
Let $J$ be a closed ideal in $L^1(K)$ and $\{e_n\}$ 
a bounded  approximate identity for $L^1(K)$,  as in Lemma \ref{theorem.1}.  
By  Corollary
 \ref{c.4} there  exists a subset $P$ of $\widehat{K}$ such that
$J=I_P$.
  Define
\begin{equation}\notag
f_P(\alpha):= \left\{
  \begin{array}{ll}
    1  & \hbox{if $\alpha \in P$,} \\[.3cm]
    0   & \hbox{if $\alpha \not \in P$.}
  \end{array}
\right.
\end{equation}
Obviously   $f_P\cdot L^2(\widehat{K})\subset L^2(\widehat{K})$ and
$\widehat{e_n}\cdot f_P$ belongs to $L^2(\widehat{K})$. Since 
the Plancherel transform is an isometry of $L^2(K)$ onto $L^2(\widehat{K})$ , there exists
$\{h_n\}$ of functions  in $L^2(K)$ such that $\widehat{h_n}=
\widehat{e_n}\cdot f_P$. Clearly $h_n\in J=I_P$ and for each $g\in
I_P$
 we have
\begin{align}\notag
\widehat{h_n\ast g}&=\widehat{h_n}\cdot \widehat{g}\\[.3cm]
\notag
&=\widehat{e_n}\cdot f_P\cdot \widehat{g}\\[.3cm]
\notag &=\widehat{e_n}\cdot \widehat{g}, \notag
\end{align}
 which implies that $h_n\ast g= e_n\ast g$. Since
 $\{e_n\}$ is a bounded approximate identity for
 $L^1(K)$, so $\{h_n\}$ is an approximate
identity for $J=I_P$.
\end{proof}

%\begin{remark}
%\emph{By previous theorem, it is straightforward to show that approximate
%identity for maximal ideals of $L^1(K)$ is bounded.
%(I have to construct one of them, it is not so easy.
%Let $I(\alpha)$ denote the maximal ideal in $L^1(K)$. }
%\end{remark}

\section{Amenability of Discrete Hypergroup Algebras}

In \cite[Theorem 2.1]{azimifard.Math.Z} we showed that if a character $\alpha$ of a polynomial hypergroup 
vanishes at infinity, then the hypergroup algebra can not be $\alpha$-amenable. In the following theorem we generalize  this fact to
 discrete hypergroups.  

\begin{theorem}\label{theorem.2}
\emph{
Let $K$ be a discrete hypergroup and $\alpha\in \widehat{K}$.  If $\alpha\in C_0(K)$,  then $L^1(K)$ is not amenable. In particular if $\pi_K(\{\alpha\})=0$,  then  $L^1(K)$ 
is not $\alpha$-left amenable. }
\end{theorem}

\begin{proof}
Let us first assume $\alpha\in C_0(K)$ with $\pi_K(\{\alpha\})=0$ and in contrary  $L^1(K)$ is 
$\alpha$-left amenable. Then by \cite[Theorem 1.2]{azimifard.Math.Z} $I(\alpha)$ has a bounded approximate identity, say $\{e_i\}_{i\in J}$ with
  $\|e_i\|_1\leq M$ for
  some $M>0$.  Let $m_\alpha$ be  the $w^\ast$-limit of $\{e_i\}$ in $L^1(K)^{\ast\ast}$. By \cite[Lemma 2]{stei}, $\{\widehat{e_i}\}$
 converges uniformly to the identity character in $\widehat{K}$ and $m_\alpha(\alpha)=0$.  
 Since $\pi_K$ is a regular measure on $\widehat{K}$, there exists an open neighbourhood
  $ {U}_\alpha$ of $\alpha$
 with $\pi_K({U}_\alpha)<\frac{\varepsilon^2}{8M^2}$, for given $\varepsilon>0$. 
 There exists a $i_0\in J$ such that $|\widehat{e_i}(\beta)-1|<\frac{\varepsilon}{\sqrt{2}}$
  for all $\beta\in {U_\alpha}^c$
 and $i\geq i_0$. Since
 \begin{equation}\notag
 |\widehat{e_i}(\beta)-1|^2\leq |\widehat{e_i}(\beta)|^2+2|\widehat{e_i}(\beta)|+
                               1\leq \|e_i\|_1^2+2\|e_i\|_1+1\leq M^2+2M+1\leq 4M^2
 \end{equation}
 for all $\beta\in \widehat{K}$, we have 
 \begin{align}\notag
 \|\widehat{e_i}-1\|_2&=\int_{\widehat{K}}|\widehat{e_i}(\beta)-1|d\pi_K(\beta)\\\notag
 &=
 \int_{{U_\alpha}}|\widehat{e_i}(\beta)-1|d\pi_K(\beta)+
 \int_{{U}^c_\alpha}|\widehat{e_i}(\beta)-1|d\pi_K(\beta)\leq \varepsilon^2.
 \end{align}
Due to   the Plancherel theorem we have  $\|e_i-\delta_e\|\rightarrow 0$ when  $i\rightarrow \infty.$ Hence 
 for every $f\in C_c(K)$  
 \begin{align}\notag
 \left|\int_K f(x)e_i(x)dm(x)-\int_K f(x)\delta_e(x)dm(x)\right|&=\left|\int_K(e_i-\delta_e)(x)f(x)dm(x)\right| \\ \notag
 &\leq \|e_i-\delta_e\|_2\|f\|_2\rightarrow 0 \hspace{1.cm}(\mbox{ as }i\rightarrow\infty).
 \end{align}
The latter 
 inequality shows that  $m_\alpha(f)=f(e)$ for all $f\in C_0(K)$. In particular  
 $m_\alpha(\alpha)=\alpha(e)=1$ which is a contradiction. Therefore $L^1(K)$ is not $\alpha$-left amenable. 
 
 Now we assume   $\pi_K(\{\alpha\})>0$.  In this 
   case $L^1(K)$ can be $\alpha$-left amenable \cite{azimifard.Math.Z}, however  we show  that $L^1(K)$ is not amenable. 
  Let ${K}':={K}\times {K}$ as
 above and  $Y:=(C_0({K}'), \|\cdot \|_\infty)$. For $f\in L^1(K)$ and $g\in Y$ define
$f\cdot g:=\pi_1f\ast g$
and 
$g\cdot f:=\pi_2f\ast g$.
It is easy to see that $Y$ is a Banach $L^1(K)$-bimodule with respect to the above 
 module multiplications.
Since $\alpha\in C_0(K)$,  $\alpha\otimes 1\in C_0(K')$  and  the maximal ideal  generated by
this character   in $M(K')$ can be regarded as a dual 
$L^1(K)$-bimodule.
 To see this,
 let   $X:=\{\varphi\in C_0(K')^\ast: \varphi(\alpha\otimes 1)=0\}$, and let  
 $\varphi\rightarrow \mu_{\varphi}$ denote   the Riesz's duality ($C_0(K')^\ast \cong  M(K')$). 
% By Riesz's theorem \cite{HewRos79I}   $C_0(K')^\ast \cong  M(K')$ 
% where the duality is given by
% \[\varphi(g)=\int_{K'} \overline{g(t_1,t_2)}
%d \mu_\varphi(t_1, t_2),\hspace{.5cm} \varphi\in C_0(K')^\ast, g\in 
% C_0(K') \mbox{ and } \mu_\varphi\in M(K'). \]
We note that since  $K'$ is discrete, 
the algebra  $L^1( K')$ can be identified with $M(K')$ via the map $f\rightarrow f m$. So, the space 
$X$ is an $L^1(K)$-submodule of $ C_0({K}')^\ast$, since
for any  $\varphi\in X$ and $f\in L^1(K)$ we have
$$
f\cdot \varphi ( \alpha\otimes 1)=\pi_2f \ast \mu_\varphi ( {\alpha\otimes 1})=
\widehat{f}(1)\widehat{\mu_\varphi}(\alpha\otimes 1)=0, $$
and likewise 
$$
 \varphi\cdot f ( \alpha\otimes 1)=\pi_1f \ast \mu_\varphi ( {\alpha\otimes 1})=
\widehat{f}(\alpha)\widehat{\mu_\varphi}(\alpha\otimes 1)=0. $$
%
%
%and likewise $\varphi\cdot f\in X$. 
Since $X$ is a (weak-$\ast$) closed subset of $C_0({K'})^\ast$, by
 \cite[Proposition 1.3]{Bade.dales.curtes} $X$ is a dual module with respect to the
module multiplications.
We can now define the continuous linear operator $D:L^1(K)\rightarrow X$ by $D(f):=\pi_1f-\pi_2f$, where for every 
$f,g\in L^1(K)$ 
\begin{align}\notag
D(f\ast g)&=\pi_1(f\ast g)-\pi_2(f\ast g)\\ \notag
&= \pi_1f\ast \pi_1g-\pi_2f\ast \pi_2g\\ \notag
&=(\pi_1f-\pi_2f)\ast \pi_1g+\pi_2f\ast (\pi_1g-\pi_2g)\\ \notag
&=D(f)\ast \pi_1g+\pi_2f\ast D(g)\\ \notag
&=D(f)\cdot g+f\cdot D(g).
\end{align}
Therefore  $D$ is a derivation.
By assumption there exists a $\varphi\in X$ such that
$D(f)=f\cdot \varphi-\varphi\cdot f$ for all $f\in L^1(K)$. % Let $\alpha_i\rightarrow 1$, as $i\rightarrow \infty$. 
Since $\widehat{K}$ separates the points of $K$ \cite{BloHey94},   there 
exists   $f\in L^1(K)$ such $\widehat{f}(\alpha)\not=\widehat{f}(1)$, however 
 \begin{align}\notag
 \widehat{f}(\alpha)-\widehat{f}(1)&=\pi_1f({\alpha\otimes 1})-\pi_2f({\alpha\otimes 1}) \\ \notag
 &
 =Df( {\alpha\otimes 1}) =
 f\cdot \varphi( {\alpha\otimes 1})-\varphi\cdot f({\alpha\otimes 1})\\ \notag
 &=
  \left(\widehat{f}(1)-\widehat{f}(\alpha)
                   \right)\widehat{\mu_\varphi}(\alpha\otimes 1)=0 \notag
  \end{align}
which is a contradiction.  Therefore $L^1(K)$ is not amenable. 
\end{proof}

\section{Examples}

      \begin{itemize}
                \item[(i)]{\bf{ Hypergroups associated to   the center of    group algebras}}

                Let $G$ be a non-abelian compact  connected Lie  group and $K$ 
                the hypergroup of conjugacy classes of $G$. 
                The center of  $L^1(G)$ is isometrically isomorphic to $L^1(K)$; see \cite{Mosak}.
                There exists a  sequence consisting of 
               irreducible unitary representations   of
                $G$ such that their dimensions  tend to infinity. Therefore,   by Theorem \ref{main.2},  
                    $L^1(K)$ is not amenable (see also \cite[Theorem.1.7]{A.S.N}).

%%%%%%%%%%%%%%%%%%%%%%%%%%%%%%%%%%%%%%%%%%%%%%%%%%%%%%%%%%%%%%%%%%%%%%%%%%%%
 \item[(ii)]{ \bf{Compact $P_\ast$-hypergroups } }

These hypergroups are due to Dunkl and Ramirez \cite{Dun.Remirez.Donald}.
 Let $0<a\leq \frac{1}{2}$ and $H_a=\NN_0\cup \{\infty\}$ denote the one point
 compactification  of $\NN_0=\NN\cup \{0\}$.
 Let $\delta_\infty$ be the identity element of $H_a$,
$\tilde{n}=n$ for all $n\in H_a$,    and define $\delta_n\ast\delta_m=\delta_{\mbox{min}(n,m)}$ for  $m\not=n\in \NN$ and 
\[
\delta_n\ast\delta_n(l) =
  \left\{
       \begin{array}{ll}
         0, & \hbox{$l<n$;} \\
         \frac{1-2a}{1-a}, & \hbox{$l=n$;} \\
         a^k, & \hbox{$l=n+k>n$.}
       \end{array}
     \right.
\]
 The Plancherel measure of $\widehat{H_a}$ is
given by
\[
\pi(\{k\})=\left\{
       \begin{array}{ll}
         1, & \hbox{$k=0$;} \\
         \frac{1-a}{a^k}, & \hbox{$k\geq 1$.}
       \end{array}
     \right.
 \]

Since $\pi(k)\rightarrow \infty$ as $k\rightarrow \infty$, by Theorem \ref{main.2}
  $L^1(H_a)$ is not amenable. Also note that from \cite{Dun.Remirez.Donald} we have 
  $\widehat{\NN_0}\setminus\{1\}\subset L^1\cap L^2(\NN_0)$, so    by  Theorem \ref{theorem.2}  $L^1(\NN_0)$ is not amenable but 
$\alpha$-left amenable in every  $\alpha\in \widehat{\NN_0}$ (see \cite{azimifard.Math.Z, Ska92}).

%%%%%%%%%%%%%%%%%%%%%%%%%%%%%%%%%%%%%%%%%%%%%%%%%%%%%%%%%%%%%%%%

\item[(iii)]{\bf{ Dual of Jacobi polynomial hypergroups }}\\
Let $K$ be Jacobi  polynomial hypergroup $\{P_n^{(\alpha,\beta)}(x)\}_{n\in \NN_0}$ of order
 $(\alpha,\beta)$, where $\alpha\geq \beta >-1$, $\alpha+\beta+1\geq 0$; see \cite{BloHey94}. 
 %Here 
 %$P_n^{(\alpha,\beta)}(x)$ are Jacobi polynomials. 
 The Haar weights are given by 
\begin{equation}\label{eq.1}
h(0)=1,\hspace{.1in}h(n)=\frac{(2n+\alpha+\beta+1)(\alpha+\beta+1)_n(\alpha+1)_n}
{(\alpha+\beta+1)n! (\beta+1)_n}, \hspace{.1in}\mbox{for} \;n\geq 1,
\end{equation}
where $(a)_n$ is the Pochhammer-Symbol. The character space of  $\NN_0$
  can be identified with $[-1,1]$ %, i.e. $\{P_{n}^{(\alpha,\beta)}(x)\}_{x\in [-1, 1]}$, 
%which 
and has the dual hypergroup structure with the Haar measure 
\begin{equation}\notag
d\pi(x)=c_{(\alpha,\beta)}(1-x)^{\alpha}(1+x)^{\beta}\chi_{[-1,1]}(x)dx \hspace{.5cm}(c_{(\alpha,\beta)}>0)
\end{equation}
where
%$C_{(\alpha,\beta)}=\frac{1}{2^{\alpha+\beta+1}}\frac{\Gamma(\alpha+\beta+2)}{\Gamma(\alpha+1)\Gamma(\beta+1)}$
%and 
$\chi$ denotes the characteristic function.

\begin{proposition}
\emph{
Let $K$ denote the compact hypergroup $[-1,1]$. Then the algebra  $L^1(K)$
 is amenable if and only if $\alpha=\beta=-\frac{1}{2}$. }
\end{proposition}
\begin{proof}

Let  $\alpha=\beta=-\frac{1}{2}$. Then  the hypergroup $[-1,1]$
is the dual of  Chebychev polynomial hypergroup with the Plancherel weights $h(0)=1$,  $h(n)=\frac{1}{2}$, $n\geq 1$.
So   by Theorem \ref{main.2}
$L^1(K)$ is amenable; see also \cite[Theorem.1.3]{A.S.N}.
In the case of   $\alpha, \beta>-\frac{1}{2}$, 
  the Plancherel weights $h(n)$ in (\ref{eq.1})  tend to infinity when  $n\rightarrow \infty$;  
  consequently, by Theorem \ref{main.2}, $L^1(K)$ is not amenable.
\end{proof}
\end{itemize}


\begin{thebibliography}{10}



 \bibitem{azimifard.Math.Z}A. Azimifard,   $\alpha$-  Amenable hypergroups. Math. Z., Vol. 262 (3) 2009.
  DOI:10.1007/s00209-009-0550-7.

\bibitem{A.S.N}
A. Azimifard, E. Samei, and N. Spronk,  
\newblock{Amenability properties of the centres of group algebras.}   J. Funct. Anal. 256 (5) (2009), 1544--1564.



\bibitem{Bade.dales.curtes}
 W. G. Bade, P. S.   Curtis  and   H. G. Dales,
\newblock  Amenability and weak amenability for Beurling and Lipschitz algebras.
\newblock {\em \emph{Proc. London Math. Soc. }} 55 (3)(1987) 359--377.


 \bibitem{BloHey94}
W.~R. Bloom and H.~Heyer,
\newblock {\em \emph{Harmonic Analysis of Probability Measures on
  Hypergroups.}}
\newblock De Gruyter, 1994.

\bibitem{BonDun73}
F.~Bonsall and J.~Duncan.
\newblock {\em \emph{Complete Normed Algebras}}.
\newblock Springer, Berlin, 1973.

\bibitem{Dun.Remirez.Donald}
C. F. Dunkl and  D. E.  Ramirez, 
\newblock{ A family of countably compact $P\sb{\ast}$-hypergroups. }
 Trans. Amer. Math. Soc.  202  (1975), 339--356. 

 
\bibitem{HewRos79I}
E. Hewitt and K. A.  Ross,
\newblock {Abstract Harmonic Analysis.}
Vol. 1,
\newblock  Springer Verlag,
  1970.
\bibitem{Jew75}
R.~I. Jewett,
\newblock Spaces with an abstract convolution of measures.
\newblock {\em \emph{Adv. in Math.}} 18(1975), 1--101.

\bibitem{Johnson}
B.E. Johnson,
\newblock Cohomology in {Banach} algebras.
\newblock {\em \emph{Amer. Math. Soc.}} 127, 1972.


\bibitem{Mosak}
J. Liukkonen and R. Mosak, 
\newblock{ Harmonic anlysis and centers of group algebras.}
\newblock{\em\emph{Trans. Amer. Math. Soc.}} 195(1974), 147--163.

 
 

  \bibitem{Ska92}
M.~Skantharajah,
\newblock Amenable hypergroups.
\newblock {\em \emph{Illinois J. Math.}} 36(1) (1992),15--46.


\bibitem{stei} U. Stegmeir,
 \newblock Centers of group algebras. Math. Ann. 243 (1979), 11--16. 
\bibitem{stei} M. Voit,
 \newblock Compact groups having almost discrete orbit hypergroups. Montash. Math. 122(3) (1996), 230--250. 
 
\bibitem{Zhang}
Y. Zhang, Approximate identities from ideals of Segal algebras on a compact group.  J. Funct. Anal. 191(2002), 123--131.
\end{thebibliography}
            \end{document}